\documentclass[a4j,12pt]{article}
\usepackage{amsmath,amsthm,amssymb,bm}

\numberwithin{equation}{section}

\newtheorem{thm}{Theorem}[section]
\newtheorem{lem}[thm]{Lemma}
\newtheorem{prop}[thm]{Proposition}
\newtheorem{cor}[thm]{Corollary}
\newtheorem{definition}[thm]{Definition}
\newtheorem{remark}[thm]{Remark}
\newtheorem{example}[thm]{Example}

\newenvironment{df}{\begin{definition}\rm}{\end{definition}}
\newenvironment{rem}{\begin{remark}\rm}{\end{remark}}

\title{Operational extreme points and  Cuntz's canonical endomorphism}

\author{Marie Choda}

\date{}

\begin{document}
\maketitle
\centerline{Osaka Kyoiku University, Osaka, Japan} 
\centerline{marie@cc.osaka-kyoiku.ac.jp}
\begin{abstract} 
Based on the fact that the Cuntz algebra $O_n$ is generated by the operators consisting of a finite operatorional partition, 
we study the  notion of operational extreme points  (which we introduce here) 
by using several completely positive maps on $O_n$.  
As a typical example, we show that the Cuntz's canonical endomorphism $\Phi_n$ 
is an operational extreme  point in the set of completely positive maps on $O_n$ 
and that it induces a completely positive map 
which is  extreme  but not operational extreme, etc. 
\end{abstract} 

keywords: {Positive linear map, convex combination, extreme point}

{Mathematics Subject Classification 2000: 46L55; 52A20, 46L37, 46L40 
\smallskip

\section{Introduction} 

From a view point of von Neumann entropy for states of the 
algebra of $n \times n$ conplex matrices $M_n(\mathbb {C})$, 
we obtained in \cite{MC1} some characterizations for unital positive   Tr-preserving  maps of  $M_n(\mathbb {C})$, 
where Tr is the standard trace of $M_n(\mathbb {C})$. 
As one of them, we showed that 
a positive unital  Tr-preserving map $\Phi$ of $M_n(\mathbb {C})$ preserves the von Neumann entropy 
of a given state $\phi$ if and only if $\Phi^* \circ \Phi$ preserves the $\phi$, 
where $\Phi^*$ is the adjoint map of $\Phi$ with respect to the Hilbert-Schmidt inner product 
$\left\langle \ \cdot\ , \ \cdot \ \right\rangle_{{\rm Tr}}$ 
for $M_n(\mathbb {C})$ 
induced by {\rm Tr}. 

This  is  based  on the property of the entropy function $\eta$ that 
$-\eta$ is an operator convex function (cf. \cite{NS}).  
\vskip 0.3cm

By keeping this operator convexity in mind,  in this paper, we study the set of completely positive maps on  unital C$^*$-algebras. 
As a generalized  notion of extreme points in the usual sense, 
we define an {\it operational extreme point} in 
the set of completely positive maps on   unital C$^*$-algebras. 
As a special case of operational extreme points, the notion of  a {\it numerical operational extreme point} appears. 
Comparing with the case of $M_n(\mathbb {C})$, in the case of infinite dimensional C$^*$-algebras, we can see different  
plays of completely positive maps. 

First we show  that a unital $*$-homomorphism of a unital C$^*$-algebra is an operational extreme point. 
\smallskip

Our notion of {\it operational convexity} is based on  the {\it finite operational partition} introduced by 
Lindblad (\cite{L}) as we describe below. 

As one of the most typical non-elementary example for a finite operatorional partition, we pick up here 
Cuntz's family of isometries and study several unital completely positive (called UCP for short) maps on $O_n, (n \ge 2)$. 
For example, the canonical shift $\Phi_n$ of the Cuntz algebra $O_n$ 
is an operational extreme point in the set of 
UCP  maps but  not  a numerical operatorional extreme point. 

In the case of $M_n(\mathbb {C})$, if $\Phi$ is a Tr-preserving $*$-homomorphism of $M_n(\mathbb {C})$, 
then the adjoint map $\Phi^*$ of  $\Phi$ is an automorphism so that $\Phi^*$ is also an  operational extreme point.  

Back to the case  of $O_n$, the unique state  of $O_n$ with the trace-like property for the 
canonical UHF-subalgebra $F_n$ of $O_n$  (that is,  $\phi_n(ab) = \phi_n (ba), \ a \in O_n, b \in F_n$ (cf. \cite{Arc})) 
plays a key role like Tr, and 
the canonical shift $\Phi_n$  is a  $\phi_n$-preserving map on $O_n$. 
We show that  the adjoint map $\Phi_n^*$ of $\Phi_n$ with respect to the $\left\langle \ \cdot \ , \cdot \ \right\rangle_{\phi_n}$ 
is nothing else but the standard left inverse $\Psi_n$ of 
$\Phi_n: \Psi_n = \Phi^*_n, \ \Phi^*_n \circ \Phi_n = {\rm id}_{O_n}$. 
However the position of  $\Phi^*_n$ is different from the case of $M_n(\mathbb {C})$ and it is not  an  operational extreme point. 

We pick up several UCP maps on $O_n$ and show relations among  operational extreme points and the notions treated in discussions related to 
extreme points in positive maps on C$^*$-algebras in \cite{St}. 

\section{Preliminaries}
Here we summarize notations, terminologies and basic facts. 

\subsection{Finite partition} 
In order to define a convex combination, we need a probability vector $\lambda = (\lambda_1, \cdots, \lambda_n )$:  
$\lambda_i \geq 0, \quad \sum_i  \lambda_i = 1$. 
Given a finite subset $x = \{x_1, ..., x_n \}$ of  a vector space $X$, 
the vector  
$\sum_i  \lambda_i x_i$ is  a convex sum of $x$ via $\lambda$. 

Now, we call such a $\lambda$  as a {\it finite partition} of 1. 
\vskip 0.3cm

Two  generalized notions of  finite partition of 1 
are given in  the framework of the non-commutative entropy as follows: 

Let $A$ be a unital C$^*$-algebra.  

\subsubsection{Finite partition of unity} 
A finite subset $\{x_1, ..., x_k \}$ of  $A$ is called  a 
{\it finite partition of unity}  by Connes-St\o rmer \cite{CS} if 
they are nonnegative operators which satisfy that 
$\sum_{i=1}^n  x_i = 1_A$. 

We denote   by $FP(A)$ the set of all  finite partitions of unity in $A$.  

\subsubsection{Finite  operational partition of unity} 
A finite subset $\{x_1, ..., x_k \}$ of  $A$ is called  a 
{\it finite  operational partition} in $A$ of unity of size $k$ by Lindblad \cite{L}     if 
$\sum_i^k x_i^* x_i = 1_A$.  

In this note, we pick up a finite subset $ \{v_1, ..., v_k \}$ of non-zero elements in $A$ such that 
$\sum_i^k v_i v_i^* = 1_A$   and call it a {\it finite operational partition of unity of size $k$} in $A$. 

The reason to use this version here is that our main target in this note is the Cuntz algebra $O_n$. 
Let us denote  the set of all  finite operational partitions of unity of size $k$ in $A$  by  
$FOP_k(A)$: 
\begin{equation*}
FOP_k(A) = \{  \{v_1, ..., v_k \} \ | \ 0 \ne v_i \in A, \forall i, \ \sum_i^k v_i v_i^* = 1_A \}
\end{equation*}
and the set of all  finite operational partition of unity in $A$ by $FOP(A)$: 
\begin{equation*}
FOP(A) = \bigcup_{k=1}^\infty  FOP_k(A). 
\end{equation*}

The group $U(A)$ of all unitaries in $A$ is a subset of the most trivial finite operational partition of 
unity with the  size $1$, that is, $U(A) \subset FOP_1(A)$. 

A linear map $\Phi$ on a unital C$^*$-algebra $A$ is  positive iff $\Phi(a)$ is  positive for all   positive $a \in A$ 
and {\it completely positive} iff $\Phi \otimes 1_k$ is  positive for all   positive integer $k$, where 
the map $\Phi \otimes 1_k$ is the map on $A \otimes M_k(\mathbb {C})$ defined by 
$\Phi \otimes 1_k (x \otimes y) = \Phi(x) \otimes  y$ for all $x \in A$ and $y \in M_k(\mathbb {C})$.  

\vskip 0.3cm 

Given an operator $v$, ${\rm Ad}\ v$ is the map given  by ${\rm Ad}\ v (x) = vxv^*$.  
If the map $\sum_{i=1}^m {\rm Ad}\ v_i$ is  unital,  
then  
$\sum_{i=1}^m v_i v_i^* = 1$, that is, the set  $\{v_1, \cdots, v_m\}$ is a  finite operational partition of the unity.

\section{\bf Operational Convex Combination}

\subsection{\bf Operational convexity} 

\begin{df} 
Let $A$ be a unital C$^*$-algebra. 
For a $\{v_i\}_{i=1}^m \in FOP(A)$ and 
a set $\{\Phi_i \}_{i=1}^m$ of linear maps  on $A$,  
we call 
the map $\sum_{i=1}^m {\rm Ad}\ v_i \circ \Phi_i$ 
an {\it operational convex combination} of $\{\Phi_i \}_{i=1}^m$ with an {\it operational coefficients} $\{v_i\}_{i=1}^m$. 
 
We also say that a subset $S$  of linear maps  on $A$ 
is  {\it operational convex} if it is closed under all operational convex combinations. 
\end{df}

\vskip 0.3cm 

\subsubsection{Operational extreme point} 
Now let us remember the  notion of  extreme points. 

Let $S$ be a convex set of a vector space $X$. 
Then a $z \in S$ is an  extreme point in $S$ if $z$ cannot be the convex combination 
of different points in $S$, 
that is, $z \in S$ is an  extreme point in $S$ if the following holds: 

\begin{eqnarray}
\lefteqn{z = \sum_{i=1}^m \lambda_i x_i, \ \{x_i\}_{i=1}^m \subset S, \ \sum_{i=1}^m \lambda_i = 1,  \  0 < \lambda_i < 1, \ \forall i} \\
 &\quad \quad \quad \quad  \quad \Longrightarrow & x_i =  z, \ \forall i \quad  (i.e.,  \lambda_i x_i = \lambda_i z, \forall i). \nonumber
\end{eqnarray}

By replacing a finite partition $\{\lambda_i\}_{i=1}^m$ of $1$ to a finite operational partition of the unity and 
a convex set to an operational convex set $S$ of mappings, 
we define an {\it operational extreme point}  of $S$.

In \cite{MC2},  we introduced  a notion of  operational extreme points 
for linear maps on  the algebra of $n \times n$ complex matrices $M_n(\mathbb {C})$. 

Here, we generalize it as follows:

\smallskip

\begin{df} 
Let $S$ be an operational convex subset of positive linear maps of a unital C$^*$-algebra  $A$ into a unital C$^*$-subalgebra  $B$ of $B(H)$. 
We say that a map $\Phi \in S$ is an {\it operational extreme point} of $S$  
if the following holds:  
\begin{eqnarray}
\lefteqn{\Phi = \sum_{i=1}^m {\rm Ad}a_i \circ \Psi_i, \ \{\Psi_i \}_{i=1}^m \subset S, \ \{a_i \}_{i=1}^m \in FOP_m(B)} \\
 &\quad \quad \quad \quad \Longrightarrow &  {\rm Ad} a_i \circ \Psi_i =  z_i \Phi, \ \text{for some} \ z_i \in \Phi(A)', 
 \quad \forall i, \nonumber
\end{eqnarray}
Here $\Phi(A)'$ is the commutant of $\Phi(A)$, i.e, $\{z \in B(H): z y = y z, \ \forall y \in \Phi(A) \}$. 
\smallskip

In the case where we can  take $\{z_i\}_{i = 1}^m$ for $\Phi$ as positive real numbers, 
we call such a  $\Phi$  a {\it numerical operational extreme point} of $S$.  
This is the case which we discussed in \cite{MC2}. 
\end{df}

\begin{rem} 
(i) 
If an operational extreme point $\Phi$ is  unital, then 
the set $\{z_i\}_{i = 1}^m$ in the definition is a finite partition  of the unity of $\Phi(A)'$, that is, 
$\{z_i\}_{i = 1}^m$ are nonnegative operators in $\Phi(A)'$ and $\sum_i z_i$ is the unity of  $\Phi(A)'$. 

In fact, each $z_i$ is nonnegative by the relation ${\rm Ad} a_i \circ \Psi_i(1) =  z_i \Phi(1) = z_i$ and 
$\sum_i  z_i = \sum_i {\rm Ad} a_i \circ \Psi_i(1) = \Phi(1) = 1.$
\smallskip

(ii) A numerical operational extreme point is a special case 
of an operational extreme point, and  an operational extreme point of UCP maps is clearly  an extreme point of UCP maps. 
As we showed in \cite {MC2} and also we show later, an  extreme point  is  not neccesaly an  operational extreme point. 
\end{rem}

\begin{prop} 
Let $A$ be a unital \rm{C}$^*$-algebra and let $\Phi$ be a unital $*$-homomorphism of $A$ into $B(H)$. 
Assume that $\Psi$ is a completely positive map of $A$ into $B(H)$ such that $\Phi - \Psi$ is also completely positive. 
Then there is a unique $z \in \Phi(A)'$ with $0 \le z \le 1_H$ such that $\Psi(x) = z\Phi(x)$ for all $x \in A$. 
\end{prop}
\smallskip
 
{\it Proof}. By a similar method to that in \cite[Section 3.5] {St} (cf. \cite{A}), we give a proof. 
Let $(\pi, V, K)$ be the minimal Stinespring representation of $\Psi$, i.e, 
$\pi$ is a representation of $A$ on $B(K)$ for a Hilbert space $K$ and 
$V: H \to K$ is a bounded linear map with $\Psi(a) = V^* \pi(a) V$ such that $K$ is spaned by $\{\pi(A)V H \}$. 
Since $\Phi - \Psi$ is CP and $\Phi$ is a *-homomorphism, we have that 
for all finite subsets $\{a_j\}_{j=1}^m \subset A$ and $\{\xi_j\}_{j=1}^m \subset H$: 
\begin{eqnarray*}
\lefteqn{\Vert \sum_{j=1}^m \pi(a_j) V \xi_j \Vert^2 = <\sum_{i,j = 1}^m V^*\pi(a_i^*a_j) V\xi_j, \xi_i>} \\
&=& <\sum_{i,j = 1}^m \Psi(a_i^*a_j) \xi_j, \xi_i> \\
 &\leq&   <\sum_{i,j = 1}^m \Phi(a_i^*a_j) \xi_j, \xi_i> \\
&=& <\sum_{i,j = 1}^m \Phi(a_i^*) \Phi(a_j) \xi_j, \xi_i> 
= \Vert \sum_{j=1}^m \Phi(a_j)  \xi_j \Vert^2.
\end{eqnarray*} 

Again, by using that $\Phi$ is a unital *-homomorphism, we have 
a unique contraction $T : H \to K$ such that  $T\Phi(a) \xi = \pi(a) V \xi$ for all $a \in A$ and $\xi \in H$. 
By taking $a = 1$, we have $T = V$ and  $T\Phi(a)  = \pi(a) T$ for  all $a \in A$. 

Let $z = T^*T$. 
Then $0 \le z \le 1$ and $z \in \Phi(A)'$ by the following: 
$$z\Phi(a) = T^*T\Phi(a) = T^* \pi(a) T = \Phi(a) T^* T = \Phi(a) z, \quad (a \in A).$$ 

On the other hand, since $T^* \pi(a) T = V^* \pi(a) V = \Psi(a)$,  it follows that 
$\Psi(a) = z\Phi(a) $ for all $a \in A$ and the uniqueness of such a $z$ comes from the fact that 
$\Phi$ is unital. 
\qed
\smallskip

The above proof shows that $z = V^*V$ for the minimal Stinespring representation $(\pi, V, K)$ of $\Psi$. 
\smallskip

As an application, we have the following: 
\smallskip

\begin{cor} A unital $*$-homomorphism of a unital {\rm C}$^*$-algebra $A$ into $B(H)$ is an 
operational extreme point in the  operational convex hull of completely positive maps of $A$ into $B(H)$. 
\end{cor}
\smallskip

{\it Proof}. Let  $\Phi$ be a unital $*$-homomorphism of a unital C$^*$-algebra $A$ into $B(H)$. 
Assume that $\Phi$ is given as an operational convex conbination:  
$\Phi = \sum_{i=1}^m {\rm Ad}\ a_i \circ \Psi_i$ with a finite operational partition $\{a_i\}_{i=1}^m$ in $B(H)$ 
and completely positive maps $\{\Psi_i\}_{i=1}^m$ of $A$ to $B(H)$. 
Then for each $i$, the two maps ${\rm Ad}\ a_i \circ \Psi_i$ and $\Phi - {\rm Ad}\ a_i \circ \Psi_i$ are completely positive maps. 
Hence  by the above proposition, 
there exists a unique $z_i \in \Phi(A)' $ with $0 \le z_i \le 1$ which satisfies that ${\rm Ad}\ a_i \circ \Psi_i = z_i \Phi$. 
This means that $\Phi$ is an operational extreme point in the operational convex set consisting of  
operational convex combinations of completely positive maps of $A$ to $B(H)$. 
\qed
\smallskip

Now let $\phi$ be a faithful state of a unital C$^*$-algebra $A$. 
Then the state $\phi$  induces the Hilbert-Schmidt inner product for $A$ by 
$$< x, y >  = \phi(y^*x), \quad (x, y \in A).$$
For a $\phi$-preserving linear map $\Psi$ of $A$,  
the adjoint map $\Psi^*$   is given by 
$$< \Psi^* (x),  y >  = < x,  \Psi (y) > , \quad (x, y \in A).$$

Let Aut$(A, \phi)$ be 
the set of all automorphisms $\Theta$ of $A$ such that $\phi \circ \Theta = \phi$. 

By remarking that the adjoint map of a *-homomorphism is not always a *-homomorphism as we give examples in the next section, 
here we show the following:  

\begin{lem} Let $\phi$ be a faithful state of a unital C$^*$-algebra $A$. 
Then for each  $\Theta \in {\rm Aut}(A, \phi)$, the adjoint map $\Theta^*$ of  $\Theta$ with respact to $\phi$ 
is in ${\rm Aut}(A, \phi) $.  
\end{lem} 

{\it Proof}. Let $\Theta \in {\rm Aut}(A, \phi)$. 
Since $< \Theta^*(x), y > = <x , \Theta(y) > = \phi(\Theta(y^*) x) = \phi(\Theta(y^*) \Theta(\Theta^{-1}(x) )) 
= \phi(\Theta(y^* \Theta^{-1}(x)) ) = \phi(y^* \Theta^{-1}(x) ) = < \Theta^{-1}(x), y > $ for all $x, y \in A$,  
we have that $\Theta^*(x) = \Theta^{-1}(x)$ for all $x \in A$ 
so that $\Theta^*(xy) = \Theta^{-1}(xy) = \Theta^{-1}(x)\Theta^{-1}(y) = \Theta^*(x) \Theta^*(y)$ for all $x, y \in A$. 
The property that $\phi \circ \Theta^* = \phi$ comes from that 
$\phi(\Theta^*(x) ) = <\Theta^*(x), 1 > = < x , \Theta(1)> = \phi(x)$. 
\qed
\smallskip

By combining this lemma and the above Corollary, we have the following: 

\begin{prop} 
Let $\phi$ be a faithful state of a unital C$^*$-algebra $A$. 
Then for each  $\Theta \in {\rm Aut}(A, \phi)$, the adjoint map $\Theta^*$ of  $\Theta$ with respact to $\phi$ 
is an operational extreme point in the  operational convex hull of completely positive maps of $A$. 
\end{prop}
\smallskip

\subsection{ Cuntz algebras} 
The Cuntz algebra $O_n$ (\cite{Cu})  is given as the C$^*$-algebra generated by 
$n (n \ge 2)$ isometries $\{S_1, \cdots, S_n \}$ on an infinite dimensional Hilbert space $H$  such that 
$\sum_i S_i S_i^* = 1_H,$  
that is,  $O_n$ is the C$^*$-algebra generated by  an operational partition of unity with the size $n$. 

Let $W^k_n$ be  the set of $k$-tuples $\mu = (\mu_1, . . . , \mu_k)$ with $\mu_m \in \{1, . . . , n\}$, 
and  $W_n$ be the union 
$\bigcup^\infty_{k=0} W^k_n$. 
If $\mu \in  W^k_n$ then $|\mu| = k$ is the length of $\mu$. 
If $\mu = (\mu_1, . . . , \mu_k) \in W_n$, then $S_\mu$ is an isometry with range projection $P_\mu = S_\mu S_\mu^*.$ 

For a given $\beta \in W^l_n$ and $i, j$ with $1 \leq i < j \leq l$, 
we let 
$$\beta_{(i, j)} = (\beta_i, \cdots, \beta_j).$$

Denote by $F_n^k$ the C$^*$-subalgebra of $O_n$ spanned by all words of the form $S_\mu S_\nu^*, \ 
\mu, \nu \in W_n^k$, which is isomorphic to the matrix algebra $M_{n^k}(\mathbb{C})$. 

The norm closure $F_n$ of $\bigcup^\infty_{k=0} F_n^k$ is the UHF-algebra of type $n^\infty$, and 
the unique tracial state $\tau_n$ of $F_n$ is extended to 
the unique state  $\phi_n$ with the trace-like property 
for  $F_n$ that $\phi_n(ab) = \phi_n(ba), (a \in O_n, b \in F_n)$ (cf. \cite{Arc}). 

The state $\phi_n$  induces the Hilbert-Schmidt inner product for $O_n$ by 
$$ \left\langle x, y \right\rangle  \ = \ \phi_n(y^*x), \quad (x, y \in O_n),$$
and the adjoint map $\Phi^*$ of a $\phi_n$-preserving linear map $\Phi$ of $O_n$ is given by 
$$\left\langle \Phi (x),  y \right\rangle  = \left\langle x,  \Phi^* (y) \right\rangle, \quad (x, y \in O_n).$$

\subsubsection{Cuntz's canonical endomorphism} 
The Cuntz's canonical endomorphism $\Phi_n$ (\cite{Cu}) is an interesting example in unital completely positive maps of infinite 
dimensional simple C$^*$-algebras, which is given as an operational convex combination of the identity: 
$$\Phi_n (x) = \sum_i S_i x S_i^*, \quad (x \in O_n).$$ 

The  map $\Psi_n$ on $O_n$ given by the form 
$$\Psi_n(x) = \frac 1n \sum_{i=1}^n S_i^*x S_i, \ (x \in O_n)$$
is called the {\it standard left inverse}  of $\Phi_n$  because $\Psi_n \circ \Phi_n$ is the identity map on $O_n$. 
The UCP map $\Psi_n$ is also  an  operational convex combination of the identity map.  

Here, we show that 
the standard left inverse $\Psi_n$ of the Cuntz canonical endmorphism $\Phi_n$ plays a role of $\Phi_n^*$. 
\vskip 0.3cm

\begin{prop} 

{\rm (i)} The Cuntz's canonical $*$-endmorphism $\Phi_n$ preserves the state $\phi_n$, that is, 
$\phi_n \circ \Phi_n = \phi_n$. 

{\rm (ii)}  $\Psi_n$ is the adjoint map of $\Phi_n$ with respect to the state $\phi_n$. 
\end{prop}

{\it Proof}. 
(i) This is trivial and we used this fact already in \cite{MC0}. 

(ii) By using the fact in \cite{Arc} that 
$ \phi_n (x) = \lim_{m \to \infty} \Psi_n^m (x)$ for all $x \in O_n$, this is shown as follows: 
\begin{eqnarray*}
\lefteqn{\left\langle  \Phi_n(x) , y  \right\rangle = \phi_n(y^* \Phi_n(x))} \\
 &=& \lim_{m \to \infty} \Psi_n^m (y^* \Phi_n(x) ) \\
 &=& \lim_{m \to \infty} \Psi_n^{m-1} (\Psi_n (y^* \Phi_n(x)) ) \\
 &=& \lim_{m \to \infty} \Psi_n^{m-1} (\frac 1n \sum_j S_j^* (y^* \sum_i S_i x S_i^*)S_j ) \\
 &=& \lim_{m \to \infty} \Psi_n^{m-1} (\Psi_n(y^*)  x ) \\
 &=& \phi_n( \Psi(y)^* x) = \left\langle  x , \Psi_n(y)  \right\rangle, \quad \text{for all} \ x, y \in O_n. 
\end{eqnarray*}
\qed
\smallskip

\begin{rem} The two kinds of properties for positive maps on  C$^*$-algebras  are treated 
in the discussions on extreme points in \cite{St}. 

One is called a {\it Jordan homomorphism} and the other is called {\it irreducible}:  

A self-adjoint linear map $\Phi$ on a C$^*$-algebra $A$ is called a {\it Jordan homomorphism} if 
$\Phi(a^2) = \Phi(a)^2$ for all self-adjoint $a \in A$. 
A Jordan homomorphism of a unital C$^*$-algebra $A$ is an extreme point of the unit ball of positive maps 
on $A$ (\cite [Proposition  3.1.5] {St}). 

A positive map $\Phi: A \to B(H)$ is called to be {\it irreducible} if $\Phi(A)'$ is the scalar operators. 
\end{rem}
\smallskip

Here we list up several properties of Cuntz's canonical endomorphism $\Phi_n$ which are related to extremalities in CP maps of $O_n$. 
\smallskip

\begin{prop} 
{\rm i)} \ The canonical endomorphism $\Phi_n$ is not irreducible. 

{\rm ii)}  The $\Phi_n$ is an operational extreme point but not a numerical operational extreme point 
in  the set of completely positive maps on $O_n$.

{\rm iii)} \ The adjoint map $\Phi_n^*$ of  $\Phi_n$ is not an 
operational extreme point in the set of completely positive maps on $O_n$. 

Moreover, $\Phi_n^*$ is not even a Jordan homomorphism.
\end{prop}

{\it Proof}. 
i)  For each $i$,  the projection $S_iS_i^*$ is contained in $\Phi_n(O_n)'$. 
In fact, for all  $x \in O_n$,  
$S_iS_i^* \sum_j S_j x S_j^* = S_i x S_i^* = \sum_j S_j x S_j^* S_iS_i^*$. 
\smallskip
 
ii) Since $\Phi_n$ is a *-endomorphism, it is an operational extreme point
by the Corollary in the previous section. 

Assume that  $\Phi_n$ is a numerical operational extreme point 
in the  operational convex hull of completely positive maps of $O_n$ into $B(H)$. 
Then we have an $n$-tupple  $\{\lambda_i \}_{i=1, \cdots, n}$ of positive real numbers such that 
$S_i x S_i^* = \lambda_i \Phi_n(x)$ for all $x \in O_n$. 
Hence $S_i S_i^* = \lambda_i$ for each $i=1, \cdots, n$, which contradicts that 
$\{S_i S_i^*; i=1, \cdots, n\}$ are mutually orthogonal projections. 
\smallskip

iii) Assume that the adjoint map $\Phi_n^*$ is an operational extreme point. 
By remembering the fact that  $\Phi_n^* $ is the left inverse 
$\Psi_n = \frac 1n \sum_{i = 1}^n {\rm Ad} \ S_i^*$ of $\Phi_n$, we have an $n$-tupple $\{z_i\}_{i = 1}^n \subset \Psi_n(O_n)'$ 
such that $\frac 1n S_i^* x S_i = z_i \Psi_n(x)$ for all $x \in O_n$, 
which implies that $z_i = \frac 1n 1 $ for all $i=1, \cdots, n$ so that $S_i^* x S_i = \Psi_n(x)$ for all $x \in O_n$ and 
$i=1, \cdots, n$. 
It does not hold.  In fact, if $j \ne i$, then  
$0 = S_i^* S_jS_j^* S_i = \Psi_n(S_jS_j^*) = \frac 1n  \sum_k S_k^* S_jS_j^* S_k = \frac 1n {\rm 1}_{O_n}$. 

Let us pick up the projection $p_1 = S_1S_1^*$. If $\Phi_n^*$ is a 
Jordan homomorphism, then  $\Phi_n^*(p_1^2) = \Phi_n^*(p_1)^2$ must hold. 
However,  the relation that $\Phi_n^* = \Psi_n$ implies that 
$\Phi_n^*(p_1^2) = \Psi_n(p_1^2) = \Psi_n(p_1) = \frac 1 n {\rm 1}_{O_n}$ and so $\Psi_n(p_1)^2 = \frac 1 {n^2}  {\rm 1}_{O_n}$, 
which contradicts that $\Psi_n(p_1)^2 = \Psi_n(p_1) = \frac 1 n {\rm 1}_{O_n}$. 
\qed
\vskip 0.3cm

\begin{rem}
A positive map $\Phi$ on a C$^*$-algebra $A$  is  said to be {\it extremal} if the only positive maps 
$\Psi$ on $A$, such that $\Phi - \Psi$ is positive, are of the form $\lambda \Phi$ with $0 \le \lambda \le 1$.  
In the set of all positive maps on $B(H)$ for a Hilbert space $H$, the map  ${\rm Ad} \ u, (u \in B(H))$ is extremal 
\cite [Proposition 3.1.3] {St}. 
\end{rem}
\smallskip

As an  example of  a completely positive map which is not a Jordan homomorphism but an 
operational extreme point in UCP maps on $O_n$, we show the following: 

\begin{prop} Assume that $\Phi$ is the map on $O_n$ given by $\Phi (x) = S_i^* x S_i$, for some $i = 1, \cdots, n$. 
Then 

{\rm (i)} $\Phi$ is a UCP map on $O_n$, which is not a Jordan homomorphism of $O_n$.

{\rm (ii)} $\Phi$ is a numerical operational extreme point in the CP maps on $O_n$. 
\end{prop}
 
{\it Proof}.  Denote by $S$ the $S_i$. 

(i) Let us consider the self-adjoint operator $S+ S^*  \in O_n$. Then 
$$\Phi((S+ S^*)^2)= S^* (S^2 + {S^*}^2 + S S^* + 1_{O_n}) S = S^2 + 2 ({1_{O_n}}) + {S^*}^2$$ 
and $(\Phi(S+ S^*) )^2 = S^2 + {S^*}^2 + SS^* + 1_{O_n}$. 
Since $S$ is an isometry but not unitary, it implys that  $\Phi((S+ S^*)^2) \ne (\Phi(S+ S^*) )^2$.  
Hence $\Phi$ is not a Jordan homomorphism. 
\smallskip

(ii) Assume that $\Phi = \sum_{i=1}^m {\rm Ad} \ a_i \circ \Psi_i$ for some integer $m$, $\{a_i \}_{i=1}^m \in FOP_m(O_n)$ and 
CP maps $\{\Psi_i \}_{i=1}^m$ of $O_n$. 
The  $\Phi$ is extremal in the positive maps on $B(H)$ (\cite{St}) so that it is extremal in the CP maps on $O_n$. 
Hence, for each $i$, we have some $\lambda_i, (0 < \lambda_i < 1)$ such that 
${\rm Ad} a_i \circ \Psi_i = \lambda_i \Phi$. This means that $\Phi$ is a numerical operational extreme point in the CP maps on $O_n$. 
\qed
\vskip 0.3cm

At the last, we show  another role of the Cuntz's canonical shift $\Phi_n$. 
The  following map $\Phi$ composed of $\Phi_n$ and ${\rm Ad} S_1^*$ is  a UCP map on $O_n$ which is 
an extreme point but not an operational extreme point of the UCP maps. 
This is an extended version of the example in \cite{MC2} for the case of matrix algebras. 
\smallskip

\begin{prop} Let $\Phi$ be the map on $O_n$ given by $\Phi =  \Phi_n \circ {\rm Ad} \ S_1^*$. 
Then 

{\rm (i)} $\Phi^2 = \Phi$ and $\Phi_n = \Phi \circ {\rm Ad}  \ S_1$. 

{\rm (ii)} $\Phi$ is an extreme point  of the set of UCP maps on $O_n$. 

{\rm (iii)} $\Phi$ is not an operational extreme point of the UCP maps on $O_n$. 
More precisely it is not numerical operational extreme.  

\end{prop} 

{\it Proof}. {\rm (i)} These are clear by the properties of the $\{S_1, ..., S_n\}$ and the relation 
${\rm Ad} \ S_1^* \circ \Phi_n = {\rm Ad} \ S_1^* S_1$.   
\smallskip

{\rm (ii)}  The following relations hold for all $\alpha, \beta \in W_n$: 
\begin{eqnarray*}
\lefteqn{\Phi (S_\alpha S_\beta^*)} \\
  &=& \left\{
\begin{array}{ll}
 0, &\alpha_1 \ne 1 \ \text{or} \ \beta_1 \ne 1 \\
 \Phi_n( S_{\alpha_{(2, |\alpha|)}} S_{\beta_{( 2, |\beta| )}  }^*), &\quad \text{otherwise} 
\end{array} 
\right.
\\
&=& \left\{
\begin{array}{ll}
 0, &\alpha_1 \ne 1 \ \text{or} \ \beta_1 \ne 1 \\
S_\alpha S_\beta^* + \sum_{i=2}^n S_i  S_{\alpha_{(2, |\alpha|)}} S_{\beta_{( 2, |\beta| )}  }^* S_i^* , &\quad  \alpha_1 = \beta_1 = 1\\
S_\alpha S_1S_1^* + \sum_{i=2}^n S_i   S_{\alpha_{(2, |\alpha|)}}    S_1 S_i^* , &\quad  \alpha_1 = 1, |\beta| = 0\\
S_1S_1^* S_\beta^*  + \sum_{i=2}^n S_i  S_1^* S_{\beta_{(2, |\beta|)}}^*  S_i^* , &\quad  |\alpha| = 0, \beta_1 = 1.
\end{array} 
\right.
\end{eqnarray*}

Now let 
$\Phi = \lambda \Psi + \lambda' \Psi', \ (0 < \lambda < 1, \ \lambda' = 1 - \lambda, \ \Psi, \Psi' \in {\rm UCP}(O_n))$. 
Then by using the standard left inverse $\Psi_n$ of $\Phi_n$, we have that 
\begin{eqnarray*}
\lefteqn{\text{Ad} \ S_1^* = \Psi_n \circ \Phi_n \circ \text{Ad} \ S_1^* = \Psi_n \circ \Phi} \\ 
 &=&  \lambda \Psi_n \circ \Psi + \lambda' \Psi_n \circ \Psi' \\
 &=& \sum_{i=1}^n \frac \lambda n \text{Ad} \ S_i^* \circ \Psi  + \sum_{i=1}^n \frac {\lambda'} n  \text{Ad} \ S_i^* \circ \Psi' 
\end{eqnarray*}
so that $\frac \lambda n \text{Ad}S_i^* \circ \Psi = \mu_i \text{Ad} \ S_1^*$ for some $0 < \mu_i < 1$ 
because $\text{Ad} \ S_1^* $ is extremal. 
Since $\text{Ad} \ S_i^* \circ \Psi$ and $\text{Ad} \ S_1^*$ are unital, we have  that $\frac \lambda n = \mu_i$, 
which implies that 
$\text{Ad} \ S_i^* \circ \Psi = \text{Ad} \ S_1^* $ for all $i$  so that 
$$\text{Ad} \ S_i^* \circ \Psi = \text{Ad} \ S_1^* =   \text{Ad} \ S_i^* \circ \Psi' = \text{Ad} \ S_i^* \circ \Phi, \ (i = 1, \cdots, n).$$ 

Remark that 
$\Phi_n$ is an operational extreme point in the completely positive maps on $O_n$ so that an extreme point. 
Since 
$$\Phi_n = \Phi_n \circ \text{Ad} \ S_1^* \circ \text{Ad} \ S_1 = \Phi \circ \text{Ad} \ S_1 
= \lambda \Psi \circ \text{Ad} \ S_1 + \lambda' \Psi' \circ \text{Ad} \ S_1,$$ 
it implies that 
$$\Phi_n = \Psi \circ \text{Ad} \ S_1 = \Psi' \circ \text{Ad} \ S_1 \quad \text{and} \quad
\Phi =  \Psi \circ \text{Ad} \ S_1 S_1^*  = \Psi' \circ \text{Ad} \ S_1 S_1^*.$$

If $\alpha_1 \ne 1$, then 
$$0 = \Phi (S_\alpha S_\beta^* S_\beta S_\alpha^*) \ge \lambda \Psi (S_\alpha S_\beta^* S_\beta S_\alpha^*) \ge 0$$ 
so that $\Psi (S_\alpha S_\beta^* S_\beta S_\alpha^*) = 0$. 
By using the Kadison-Schwartz inequality 
$$0 = \Psi (S_\alpha S_\beta^* S_\beta S_\alpha^*) \ge \Psi (S_\alpha S_\beta^*) \Psi( S_\beta S_\alpha^*)$$ 
which implies that 
$\Psi (S_\alpha S_\beta^*) = 0$. 
Similarly, $\Psi (S_\alpha S_\beta^*) = 0$ if $\beta_1 \ne 1$ so
that 
$$\Psi (S_\alpha S_\beta^*) = 0 = \Psi' (S_\alpha S_\beta^*) \quad \text{if} \ \alpha_1 \ne 1 \ \text{or} \ \beta_1 \ne 1.$$

If $\alpha_1 = 1$, then 

$\Psi (S_\alpha S_1 S_1^*)
= \Psi (S_1 S_{(\alpha_2, |\alpha|)} S_1 S_1^*) = \Psi \circ \text{Ad}S_1(S_{(\alpha_2, |\alpha|)}S_1).$ 

On the other hand, since $\Psi \circ \text{Ad}S_1 = \Phi_n$ and $\Psi (S_\alpha S_i S_i^*) = 0$ if $i \ne 1$, 
we have that 
\begin{eqnarray*}
\lefteqn{\Psi (S_\alpha) = \Psi (S_\alpha S_1 S_1^*) + \sum_{i=2}^n \Psi (S_\alpha S_i S_i^*) = \Psi (S_\alpha S_1 S_1^*) } \\
 &=& \Psi \circ \text{Ad}S_1(S_{(\alpha_2, |\alpha|)}S_1) = \Phi_n(S_{(\alpha_2, |\alpha|)}S_1) 
= \Phi(S_\alpha)
\end{eqnarray*}
and  a similar relation holds  for $\Psi'$ so that 
$$\Psi (S_\alpha)  = \Phi(S_\alpha) =  \Psi' (S_\alpha) \quad \text{if} \ \alpha_1 =1.$$

If $\alpha_1 = 1 = \beta_1$, then $S_\alpha S_\beta^* = S_1 S_1^*(S_\alpha S_\beta^* )S_1 S_1^*$. 
Hence by the relation that $\Psi \circ \text{Ad}S_1 S_1^* = \Psi' \circ \text{Ad}S_1 S_1^* = \Phi$, we have that 
$$\Psi (S_\alpha S_\beta^* ) =  \Phi (S_\alpha S_\beta^* )  =  \Psi' (S_\alpha S_\beta^* ) \quad \text{if} \ \alpha_1 = \beta_1 = 1.$$

As a consequence, these relations show that the map $\Phi$ is an extreme point  of the set of UCP maps on $O_n$. 
\smallskip

(iii) 
The $\Phi$ is given as an operational convex combination  $\sum_{i=1}^n \text{Ad}S_i \circ \text{Ad}S_1^*$ 
via a finite operational partition $\{S_1, \cdots, S_n\}$ and the UCP map $\text{Ad}S_1^*$. 
If $\Phi$ is  a numerical operational extreme point of the UCP maps on $O_2$, then   
there exists $\{\lambda_i \}_{i=1}^n$ such that $0 < \lambda_i < 1$, 
$\sum_{i=1}^n \lambda_i = 1$ and $\text{Ad}S_i S_i^* = \lambda_i \Phi$ for all $i$.  
Hence $S_iS_i^* =  \lambda_i 1_{O_n}$ for all $i$ which 
contradict the definition of $\{S_1, \cdots, S_n\}$. 
\qed

\end{document}